\documentclass{amsart}

\newtheorem{theorem}{Theorem}[section]
\newtheorem*{maintheorem}{Theorem}
\newtheorem{lemma}[theorem]{Lemma}

\newtheorem{corollary}[theorem]{Corollary}

\theoremstyle{definition}

\newtheorem*{acknowledgement}{Acknowledgement}

\theoremstyle{remark}

\DeclareFontFamily{U}{wncy}{}
\DeclareFontShape{U}{wncy}{m}{n}{<->wncyr10}{}
\DeclareSymbolFont{mcy}{U}{wncy}{m}{n}
\DeclareMathSymbol{\Sh}{\mathord}{mcy}{"58}

\newcommand\mynote[1]{\marginpar{\ \\ \small \tt #1}}
\newcommand\bel[1]{{\mynote{#1}}\begin{equation}\label{#1}}
\newcommand\mylabel[1]{\label{#1}}

\newcommand  {\shL}     {\mathcal{L}}

\newcommand  {\Fr}      {\operatorname{Fr}}

\renewcommand{\O}       {\mathcal{O}}

\newcommand  {\ra}      {\rightarrow}

\newcommand  {\Spec}    {\operatorname{Spec}}

\def\mydate{\number\day\space\ifcase\month \or January\or February\or March\or 
April\or May\or June\or July\or
August\or September\or October\or November\or December\fi \space\number\year}

\DeclareFontFamily{U}{wncy}{}
\DeclareFontShape{U}{wncy}{m}{n}{<->wncyr10}{}
\DeclareSymbolFont{mcy}{U}{wncy}{m}{n}
\DeclareMathSymbol{\Sh}{\mathord}{mcy}{"58}

\usepackage{amscd,amssymb,amsmath}

\begin{document}

\title[Genus-change]
      {On genus-change in algebraic curves over nonperfect fields}

\author[Stefan Schr\"oer]{Stefan Schr\"oer}
\address{Mathematisches Institut, Heinrich-Heine-Universit\"at,
40225 D\"usseldorf, Germany}
\curraddr{}
\email{schroeer@math.uni-duesseldorf.de}


\dedicatory{5 March 2007}

\begin{abstract}
I give a new proof, in scheme-theoretic language,
of Tate's old result on genus-change over
nonperfect fields in characteristic $p>0$.
Namely, for normal geometrically integral   curves,
the difference between   arithmetic and geometric genus
over the algebraic closure is divisible by $(p-1)/2$.
\end{abstract}

\maketitle

\section*{Introduction}
\mylabel{introduction}

A distinctive feature of geometry in characteristic $p>0$
is that  a regular scheme  $X$ of finite type over a nonperfect field $K$
may cease to be regular after   purely inseparable base change.
This striking behavior   easily appears for the generic fiber of     morphisms
$f:S\ra B$ between smooth schemes  over   algebraically closed ground fields:
Here $K=\kappa(B)$ is the function field of $B$, and $X=S_K$ is the
generic fiber in the sense of scheme theory. 
When it comes to classification of  fibrations,
for example in  the Enriques classification of surfaces, the theory of Albanese maps,
or the minimal model program, it is   crucial to understand this behavior.

The simplest situation is that $X$ is a proper normal curve over a nonperfect field
$K$. If $K\subset K'$ is a purely inseparable field extension,
the induced curve $X'=X\otimes_K K'$ is not necessarily normal.
Let $\tilde{X}'$ be its normalization.
Then the genus $\tilde{g}=h^1(\O_{\tilde{X}'})$ may be strictly smaller that the
genus $g=h^1(\O_X)$ of our original curve.
Tate \cite{Tate 1952} proved that such genus-change   is not arbitrary:

\begin{maintheorem}
(Tate) The difference $g-\tilde{g}$ is divisible by $(p-1)/2$.
\end{maintheorem}

In particular, this puts an upper bound on the characteristic
in terms of  the  possible genera  occurring in genus-change situations.
This   is a prominent manifestation of the intuitive principle
that a given geometrical deviation in positive characteristics in a fixed dimension
should occur only at finitely many primes. Example: Quasielliptic fibrations
(the case $g=1$, $\tilde{g}=0$) are possible only at prime $p=2$ and $p=3$.

Back in 1952,
Tate naturally stated and proved his   result in the language  of function fields and repartitions. 
In my opinion, it is desirable to have a proof
in the modern  language of schemes as well.
In the special case $\tilde{g}=0$,
Shepherd-Barron \cite{Shepherd-Barron 1991} found such a proof
for the inequality $g\geq (p-1)/2$,
using vector bundles on algebraic surfaces.
The goal of this paper is to give an easy direct proof of Tate's result,
using     relative dualizing sheaves and
relative Frobenius maps for curves. The result essentially takes the following
form:

\begin{maintheorem}
Let $Y$ be the normalization of the Frobenius pullback $X^{(p)}$.
Then the degree of the relative dualizing sheaf
$\omega_{Y/X^{(p)}}$ is divisible by $p-1$.
\end{maintheorem}

Our proof hinges on a result of Kiehl and Kunz \cite{Kiehl; Kunz 1965}, which implies that
a finite universal homeomorphism between regular curves admits
locally $p$-bases.
I expect that this approach  should  yield  result in higher dimensions as well.
The paper also contains some results on normalization of geometrically integral schemes
after Frobenius pullbacks.

\begin{acknowledgement}
I wish to thank Igor Dolgachev for stimulating discussions.
\end{acknowledgement}

\section{Normalization after Frobenius pullback}
\mylabel{normalization}

Let $K$ be a field of characteristic $p>0$,
and $X$ be a   normal $K$-scheme of finite type.
Throughout, we assume that $X$ geometrically integral,
that is, the induced schemes $X'=X\otimes_K K'$ remain  integral for all base field extensions $K\subset K'$.
The scheme $X'$, however, is not necessarily normal.
In this section, we shall collect some useful facts about the normalization of $X'$.
Our first observation is:

\begin{lemma}
\mylabel{geometrically integral}
Let $K\subset K'$ be a base field extension.
Set $X'=X\otimes_K K'$, and let $Y'\ra X'$ be the normalization.
Then the $K'$-scheme $Y'$ is geometrically integral.
\end{lemma}

\proof
Geometric irreducibility  and geometric reducedness easily follow  from \cite{EGA IVb}, Proposition 4.5.9,
and Proposition 4.6.1, respectively.
\qed

\medskip
Next, we consider   the   Frobenius pullback
$X^{(p)}$, which is defined by the cartesian square
$$
\begin{CD}
X^{(p)} @>>> X\\
@VVV @VVV\\
\Spec(K) @>>F_K> \Spec(K).
\end{CD}
$$
Here $F_K$ denotes the  absolute Frobenius morphism, which corresponds to  the Frobenius map
$\Fr:K\ra K$, $\lambda\mapsto\lambda^p$. 
The $K$-scheme  $X^{(p)}$ is  of finite type and geometrically integral,
but not necessarily normal.
In any case, the Frobenius pullback is closely related to the original   normal scheme $X$ via
the relative Frobenius morphisms $F_{X/K}:X\ra X^{(p)}$. 
This   is a finite  universal homeomorphism,
coming from  the commutative square
$$
\begin{CD}
X @>F_X>> X\\
@VVV @VVV\\
\Spec(K) @>>F_K> \Spec(K).
\end{CD}
$$
Using   iterated Frobenius maps, we obtain 
similarly the iterated Frobenius pullback $X^{(p^n)}$,
together with the iterated relative Frobenius morphism $F^n_{X/K}:X\ra X^{(p^n)}$.
The following observation will be useful:

\begin{lemma}
\mylabel{geometrically normal}
There is an integer $n_0\geq 0$ such that
for all integers $n\geq n_0$  the normalization
of $X^{(p^n)}$ is geometrically normal.
\end{lemma}

\proof
Clearly, it suffices to find one integer $n\geq 0$
so that the normalization of $X^{(p^n)}$ is geometrically normal.
To do so, choose a perfect closure $K\subset K^{p^{-\infty}}$.
Set $Z=X\otimes_K K^{p^{-\infty}}$, and let $\nu:\tilde{Z}\ra Z$ be the normalization.
The  scheme $\tilde{Z}$ is geometrically normal over $K^{p^{-\infty}}$, because
the latter   is perfect.
According to \cite{EGA IVc}, Theorem 8.8.2,  there is 
an intermediate field $K\subset K'\subset K^{p^{-\infty}}$ that is finite over $K$,
so that the scheme $\tilde{Z}$ and the morphism $\nu:\tilde{Z}\ra Z$ over $K^{p^{-\infty}}$ 
are induced
from a  scheme $\tilde{X}'$ and a morphism $\nu':\tilde{X}'\ra X'$ over $K'$.
Here of course we write $X'=X\otimes_K K'$.
By \cite{EGA IVb}, Corollary 6.7.8, the $K'$-scheme $\tilde{X}'$ is geometrically normal.
Since $\nu'$ is birational, $\nu'$   must be the 
normalization map of $X'$, and remains so
after any base field extension of $K'$.

Since the field extension $K\subset K'$ is finite and purely inseparable,
there is an integer $n\geq 0$ with the property $\lambda^{p^n}\in K$
for all $\lambda\in K'$. By the universal property of splittings fields, there
exists a homomorphism $i:K'\ra K$ so that the composite
$K\subset K' \stackrel{i}{\hookrightarrow} K$
equals the $n$-fold Frobenius map.
Consequently, $\tilde{X}'\otimes_{K'} K$ 
 is the normalization of $X^{(p^n)}=X'\otimes_{K'} K$,
where the tensor products are with respect to $i:K'\ra K$.
This concludes the proof, since we saw in the preceding paragraph that  $\tilde{X}'$ is geometrically normal.
\qed

\section{Genus-change for algebraic curves}
\mylabel{genus-change}

Now let $X$ be a proper normal curve over $K$.
As in the preceding section, we assume that $X$  is geometrically integral.
The degree of an invertible sheaf $\shL$ on $X$ is defined
as the integer $\deg(\shL)=\chi(\shL)-\chi(\O_X)$.
The main result of this paper relates the degrees of the
dualizing sheaves on the Frobenius pullback  and its normalization.
I formulate it in terms of  the relative dualizing sheaf:

\begin{theorem} 
\mylabel{degree divisible}
Let $\nu:Y\ra X^{(p)}$ be the normalization map.
Then the degree of the relative dualizing sheaf
$\omega_{Y/X^{(p)}}$ is divisible by $p-1$.
\end{theorem}

\proof
The idea is to compute with relative dualizing sheaves on  $X$.
Since $X$ is normal, there is a unique morphism
$f:X\ra Y$ with $F_{X/K}=\nu\circ f$. The various relative dualizing sheaves
satisfy
\begin{equation}
\label{tensor product}
\omega_{X/X^{(p)}} = \omega_{X/Y}\otimes f^*(\omega_{Y/X^{(p)}}).
\end{equation}
Similarly we have
$\omega_{X/K}=\omega_{X/X^{(p)}} \otimes F^*_{X/K}(\omega_{X^{(p)}/K})$.
Together with the formula  $F^*_{X/K}(\omega_{X^{(p)}/K})=F_X^*(\omega_{X/K})=\omega_{X/K}^{\otimes p}$, this yields
\begin{equation}
\label{first}
\omega_{X/X^{(p)}}=\omega_{X/K}^{\otimes (1-p)}.
\end{equation}
On the other hand,  Kiehl and Kunz  proved that the morphism
$f:X\ra Y$ admits locally $p$-bases (\cite{Kiehl; Kunz 1965}, Korollar 2 of Satz 5). Therefore 
the sheaf of relative K\"ahler differentials $\Omega^1_{X/Y}$ is locally free
of finite rank. It is related to the relative dualizing sheaf by 
\begin{equation}
\label{second}
\omega_{X/Y}=\det(\Omega^1_{X/Y}) ^{\otimes (1-p)},
\end{equation}
according to loc.\ cit., Satz 9.
Substituting formula (\ref{second}) and (\ref{first}) into   (\ref{tensor product}),
we infer that the degree of $ f^*(\omega_{Y/X^{(p)}})$ is divisible by $p-1$.
Finally, observe that  
$
\deg(f^*(\omega_{Y/X^{(p)}})) =\deg(f)\cdot\deg(\omega_{Y/X^{(p)}}), 
$
by the projection formula. Clearly, the surjection $f:X\ra Y$ is   purely inseparable, hence
its degree is a $p$-power.
From this we infer that  the degree of $\omega_{Y/X^{(p)}}$ must be divisible by $p-1$.
\qed

\medskip
Actually, the preceding result is equivalent to the following  seemingly stronger
statement:

\begin{theorem}
\mylabel{arbitrary extension}
Let $K\subset L$ be an arbitrary field extension, and $Y$ be the
normalization of $X\otimes_K L$. Then the degree of the relative dualizing sheaf $\omega_{Y/X\otimes_K L}$ is divisible by $p-1$.
\end{theorem}

\proof
First note   the following transitivity property:
Suppose $K\subset L'\subset L$ is an intermediate field. Let $Y'$
be the normalization of $X\otimes_K L'$. Then $Y$ is the normalization
of both $X\otimes_K L$ and $Y'\otimes_{L'} L$, and we have
$$
\omega_{Y/X\otimes_K L} =
\omega_{Y/Y'\otimes_{L'}L} \otimes \varphi^*(\omega_{Y'/X\otimes_K L'}\otimes_{L'}L),
$$
where $\varphi:Y\ra Y'\otimes_{L'} L$ is the normalization map.
Clearly, if two of the three dualizing sheaves have degree divisible
by $p-1$, so has the third.

Using this transitivity property, we now settle the special case
that the field $L$ in our field extension $K\subset L$ is perfect.
Choose an integer $n\geq 0$ so that the normalization $Y'$
of the Frobenius pullback $X^{(p^n)}$ is geometrically normal,
as in Lemma \ref{geometrically normal}. Since the field $L$ is perfect,
there exists precisely one homomorphism $i:K\ra L$ so that the 
composition $i\circ \Fr^n$ is our given extension $K\subset L$, according
to \cite{A 4-7}, Chap.\ V, \S 5, No.\ 2, Proposition 3.
Consider the intermediate field $L'=i(F)$.
Induction on $n$, together with the    transitivity property and
Theorem \ref{degree divisible}, shows that the degree of $\omega_{Y'/X\otimes_K L'}$ is divisible by $p-1$.
Since $Y'$ is geometrically normal, we have $Y=Y'\otimes_{L'}L$.
Another application of the transitivity property yields that
the degree of $\omega_{Y/X\otimes_K L}$ is divisible by $p-1$.

It remains to treat   the general case. Choose a perfect closure
$L\subset L^{p^{-\infty}}$. According to the preceding
paragraph, the theorem holds true for the   field extensions
$K\subset L^{p^{-\infty}}$ and $L\subset L^{p^{-\infty}}$.
By   transitivity, it  must hold   for $K\subset L$ as well.
\qed

\medskip
We now may retrieve Tate's result:

\begin{corollary}
\mylabel{tate's result}
(Tate) 
Let $K\subset L$ be an arbitrary field extension,
and $Y$ be the normalization of $X\otimes_K L$.
Then  the difference  $h^1(\O_X)-h^1(\O_Y)$ is divisible by $(p-1)/2$.
\end{corollary}

\proof
 According to Lemma \ref{geometrically integral},
the $L$-scheme $Y$ is geometrically integral.
In particular, we have 
$K=H^0(X,\O_X)$ and $L=H^0(Y,\O_Y)$.
Whence
$$
h^1(\O_X)-h^1(\O_Y)=\chi(\O_Y)-\chi(\O_X) =
\frac{1}{2}(\deg(\omega_{X/K}) -\deg(\omega_{Y/K})),
$$
the latter by Serre duality and Riemann--Roch. 
The term on the right is nothing but $-\frac{1}{2}\deg(\omega_{Y/X\otimes_K L})$,
so the statement follows from Theorem \ref{arbitrary extension}.
\qed


\end{document}